\documentclass[11pt,twoside]{amsart}

\usepackage[dvipsnames]{xcolor}
\usepackage{hyperref}
\hypersetup{colorlinks=true, linkcolor=Blue, citecolor=Maroon}
\usepackage{amsthm, amsmath, amscd, amssymb,centernot,txfonts}
\usepackage{tikz-cd}
\usepackage{tikz}

\usepackage{lipsum}
\usepackage{multicol}
\usepackage{amsfonts}
\usepackage{adjustbox}
\usepackage{mathrsfs}
\usepackage[all]{xy}
\usepackage[left=2cm,top=2.5cm,bottom=3cm,right=3cm]{geometry}
\usepackage{dirtytalk}
\usepackage{mathtools}
\usepackage{enumitem}
\usepackage{comment}
\usepackage{fourier} 
\usetikzlibrary{shapes,arrows}

\theoremstyle{plain}
\numberwithin{equation}{section}
\newtheorem{theorem}{Theorem}[section]
\newtheorem{proposition}[theorem]{Proposition}
\newtheorem{lemma}[theorem]{Lemma}
\newtheorem{corollary}[theorem]{Corollary}

\newtheorem{definition}[theorem]{Definition}

\newtheorem{set-up}[theorem]{Set-up}

\theoremstyle{definition}
\newtheorem{remark}[theorem]{Remark}

\newcommand*{\QEDB}{\hfill\ensuremath{\square}}

\tikzstyle{decision} = [diamond, draw, , 
    text width=4.5em, text badly centered, node distance=3cm, inner sep=0pt]
\tikzstyle{block} = [rectangle, draw, , 
    text width=10em, text centered, rounded corners, minimum height=2em]
\tikzstyle{block1} = [rectangle, draw, , 
    text width=5em, text centered, rounded corners, minimum height=2em]    
\tikzstyle{line} = [draw, -latex']
\tikzstyle{cloud} = [draw, ellipse,, node distance=3cm,
    minimum height=2em]
\definecolor{light-gray}{gray}{0.95}    
\begin{document}

\title[K3 carpets on minimal rational surfaces and their smoothings]{K3 carpets on minimal rational surfaces and their smoothings}

\author[P. Bangere]{Purnaprajna Bangere}
\address{Department of Mathematics, University of Kansas, Lawrence, USA}
\email{purna@ku.edu}

\author[J. Mukherjee]{Jayan Mukherjee}
\address{Department of Mathematics, University of Kansas, Lawrence, USA}
\email{j899m889@ku.edu}

\author[D. Raychaudhury]{Debaditya Raychaudhury}
\address{Department of Mathematics, University of Kansas, Lawrence, USA}
\email{debaditya@ku.edu}


\maketitle
\begin{abstract}

In this article, we study K3 double structures on minimal rational surfaces $Y$. The results show there are infinitely many non-split abstract K3 double structures on $Y = \mathbb{F}_e $ parametrized by $\mathbb P^1$, countably many of which are projective. For $Y = \mathbb{P}^2$ there exists a unique non-split abstract K3 double structure which is non-projective (see \cite{Dr}). We show that all projective K3 carpets can be smoothed to a smooth K3 surface. One of the byproducts of the proof shows that unless $Y$ is embedded as a variety of minimal degree, there are infinitely many embedded K3 carpet structures on $Y$. Moreover, we show any embedded projective K3 carpet on $\mathbb F_e$ with $e<3$ arises as a flat limit of embeddings degenerating to $2:1$ morphism. The rest do not, but we still prove the smoothing result. We further show that the Hilbert points corresponding to the projective K3 carpets supported on $\mathbb{F}_e$, embedded by a complete linear series are smooth points if and only if $0\leq e\leq 2$. In contrast, Hilbert points corresponding to projective (split) K3 carpets supported on $\mathbb{P}^2$ and embedded by a complete linear series are always smooth. The results in \cite{BGG} show that there are no higher dimensional analogues of the results in this article. 
\end{abstract}

\section{Introduction} A K3 carpet on a regular surface $Y$ is a locally Cohen-Macaulay double structure on $Y$ with the same invariants as a smooth K3 surface (i.e., regular with trivial canonical sheaf).  In this article we study the relationship between deformation theory of double covers and K3 carpets on Hirzebruch surfaces and $\mathbb P^2$.  A special case of this relationship concerning hyperelliptic K3 surfaces and double structures on rational normal scrolls was studied in \cite{GP97}. Multiple structures arise in a variety of contexts in algebraic geometry. For example, they arise in the study of vector bundles and linear series (see \cite{BM}, \cite{HV}, \cite{M1}), as well as in deformation and moduli problems (see \cite{BGG}, \cite{F}, \cite{GGP2}, \cite{GGP08}). One of the interesting features of K3 carpets is that the hyperplane section of K3 carpets supported on scrolls are canonical ribbons. The study of canonical ribbons was first proposed by Bayer and Eisenbud (see \cite{BE}).  
It follows from the works of Eisenbud and Green (see \cite{EG}) that an affirmative answer to Green's conjecture in the case of canonical ribbons will imply Green's conjecture for general curves. The reader is directed to \cite{D}, \cite{Dr}, \cite{RS}, and \cite{V} for some further reading in this direction.

\smallskip

In this article we show that a K3 double structure on a minimal rational surface $Y$ (and $\mathbb F_1$) is in general not projective, unlike ribbons on curves.   
The situation in this article resembles more that of K3 double structures arising from unramified K3 double covers of Enriques surfaces than ramified covers of Hirzebruch surfaces embedded as a variety of minimal degree, even though covers that appear in this article are all ramified covers. We show that the non-split K3 carpets on a Hirzebruch surface $Y$ are parametrized by the projective line $\mathbb{P}^1$ and the locus that parametrizes the projective K3 carpets is an infinite countable set. This invokes the notion of projective and non-projective K3 surfaces, where the former lie on infinite countably many codimension one families in the moduli space of K3 surfaces. In \cite{GP97} it was shown that there is a unique K3 carpet supported on rational normal scrolls (see \cite{GP97}, Proposition 1.7) and that it can be realized as the limit of  embedded models of smooth polarized K3 surfaces. The homogeneous ideals of these K3 carpets have been described by Eisenbud and Schreyer in \cite{ES}. The following result shows the sharp contrast between embedded K3 carpets on minimal and non-minimal degree embeddings of $Y$. 

\begin{theorem}\label{it2}
(See Theorems ~\ref{counting carpets} and ~\ref{countingp2})
\begin{itemize}[leftmargin=0.3in]
\item[$(1)$]Suppose $j: S \hookrightarrow \mathbb{P}^n$ be an embedding of a Hirzebruch surface $S$ induced by the complete linear series of $\mathcal{O}_S(1) = \mathcal{O}_S(aC_0 + bf)$ (and hence $b \geq ae + 1$ and $n+1 = h^0(\mathcal{O}_S(1))$). Let $N+1 \geq h^0(\mathcal{O}_S(1))$ and let $k: \mathbb{P}^n \hookrightarrow \mathbb{P}^N$ be a linear embedding. Let $i = k \circ j$. 
Then the K3 carpets on $S$ with an embedding inside $\mathbb{P}^N$ extending the embedding $i$ are parametrized by the non-empty Zariski open subset $U$ of $\,\,\mathbb{P}(H^0(\mathcal{N}_{S/\mathbb{P}^N} \otimes K_S))$, where $\mathbb{P}(H^0(\mathcal{N}_{S/\mathbb{P}^N} \otimes K_S))$ is the projective space of lines in $H^0(\mathcal{N}_{S/\mathbb{P}^N} \otimes K_S)$ and the points of $U$ correspond to nowhere vanishing global sections of $\mathcal{N}_{S/\mathbb{P}^N} \otimes K_S$. The dimension of $H^0(\mathcal{N}_{S/\mathbb{P}^N} \otimes K_S)$ is given by $$h^0(\mathcal{N}_{S/\mathbb{P}^N} \otimes K_S) = (N+1)\frac{(a-1)}{2}(2b-2-ae) + 1.$$ 
\item[$(2)$] Suppose that $i: S = \mathbb{P}^2 \hookrightarrow \mathbb{P}^n$ be an embedding induced by the complete linear series $\mathcal{O}_{\mathbb{P}^2}(d)$ (and hence $n+1 = h^0(\mathcal{O}_{\mathbb{P}^2}(d))$). Let $N+1 \geq h^0(\mathcal{O}_{\mathbb{P}^2}(d))$ and let $k: \mathbb{P}^n \hookrightarrow \mathbb{P}^N$ be a linear embedding. Let $i = k \circ j$. 
Then the K3 carpets on $S$ with an embedding inside $\mathbb{P}^N$ that extends the embedding $i$ are parametrized by the non-empty Zariski open subset $U$ of $\,\,\mathbb{P}(H^0(\mathcal{N}_{S/\mathbb{P}^N} \otimes K_S))$, where $\,\,\mathbb{P}(H^0(\mathcal{N}_{S/\mathbb{P}^N} \otimes K_S))$ is the projective space of lines in $H^0(\mathcal{N}_{S/\mathbb{P}^N} \otimes K_S)$ and the points of $U$ correspond to nowhere vanishing global sections of $\mathcal{N}_{S/\mathbb{P}^N} \otimes K_S$. The dimension of $H^0(\mathcal{N}_{S/\mathbb{P}^N} \otimes K_S)$ is given by  $$h^0(\mathcal{N}_{S/\mathbb{P}^N} \otimes K_S) = \frac{1}{2}(N+1)(d-1)(d-2).$$
\end{itemize}
\end{theorem}
Note that if in the above theorem we assume that $i: S = \mathbb{F}_e \hookrightarrow \mathbb{P}^N $ is an 
embedding induced by the complete linear series of $\mathcal{O}_S(aC_0+bf)$, we have that the K3 carpets embedded in $\mathbb{P}^N$ supported on $S$ are parametrized by a non-empty open set inside a projective space of dimension $$\frac{1}{4}(a^2-1)((2b-ae)^2-4).$$ 
Note that the result of \cite{GP97} is a special case of the above theorem. In \cite{GP97}, the authors show that there is an unique embedded K3 carpet where $S$ is embedded as a surface of minimal degree: recall that $S$ is embedded as a minimal degree surface if and only if $a=1$ (see \cite{EH}), in which case the formula above yields 0.

The above theorem brings up new issues unseen for multiple structures on surfaces in general and on minimal degree embeddings of $\mathbb{F}_e$ in particular. First is the existence of abstract K3 double structures on minimal rational surfaces $Y$ (and $\mathbb{F}_1$) and existence of projective K3 double structures among them. The second is the number of embedded K3 double structures on an arbitrary embedding of $Y$ and their smoothing.  The existence/non-existence of embedded non-split multiple structures on a given embedded variety and their smoothing is closely related to the deformation theory of finite covers (see \cite{GGP2}). In Section 3, Theorem ~\ref{smoothing1} shows that for every abstract projective K3 carpet on $Y=\mathbb{F}_e$ with $0\leq e\leq 2$, there is an associated double cover that is  a K3 surface.  We show that this double cover can be deformed to an embedding, thereby smoothing the K3 carpet on $Y$. For $e>2$, there are no such K3 double covers even though there are K3 double structures on the corresponding Hirzebruch surfaces $\mathbb{F}_e$, but we show these double structures are smoothable, as well.  

\begin{theorem}\label{B}(Theorems ~\ref{smoothing1}, ~\ref{smoothing2}, ~\ref{smoothingp2}, ~\ref{sm1} and ~\ref{sm2})
Let $\widetilde{S}$ be a projective $K3$ carpet supported on a minimal rational surface $S$. Then $\widetilde{S}$ is smoothable. Moreover, the following results hold.
\begin{itemize}
    \item[$(1)$] Let $\widetilde{S}$ be a projective $K3$ carpet supported on a Hirzebruch surface $S = \mathbb{F}_e$ and embedded inside $\mathbb{P}^N$ by the complete linear series of a very ample line bundle. Then $\widetilde{S}$ is a smooth point of its Hilbert scheme if and only if $0\leq e\leq 2$.
    \item[$(2)$] Let $\widetilde{S}$ be a projective $K3$ carpet supported on $S = \mathbb{P}^2$ and embedded inside $\mathbb{P}^N$ by the complete linear series of a very ample line bundle. Then $\widetilde{S}$ is a smooth point of its Hilbert scheme.
\end{itemize}
\end{theorem}

In a recent article \cite{BGG}, it is shown that there are no higher dimensional analogues of these results. In  \cite{BGG} authors introduce the notion of generalized hyperelliptic  varieties. It has been shown that the deformations of \textit{generalized hyperelliptic polarized Calabi-Yau varieties} of dimension $\geq 3$ are again generalized hyperelliptic varieties (\cite{BGG}, Theorem 4.5). They also show there are no Calabi-Yau double structures on higher dimensional scrolls, and more generally on higher dimensional projective bundles. However the results in this article show that a general deformation of a generalized hyperelliptic K3 surface is no longer generalized hyperelliptic. \textcolor{black}{It is also interesting to compare our situation with generalized hyperelliptic polarized Fano-K3 pairs. Notice that \cite{BGG}, Definition 2.8 only defines a Fano-K3 pair of index $m\geq 3$. However, if we put $m=2$ in the definition, then it yields a generalized hyperelliptic polarized K3 surface. It follows from \cite{BGG}, Corollary 2.3, that for a Fano-K3 pair $(X,L)$ of dimension $m\geq 3$, there is no non-split ribbon with conormal bundle same as the trace zero module of the morphism $\varphi$ induced by the polarization, and supported on the image of $\varphi$ unless the image is a hyperquadric in $\mathbb{P}^{m+1}$.} 
 
Theorem ~\ref{B} (1) generalizes \cite{GP97}, Theorem 4.1. Since it shows that a projective K3 carpet $\widetilde{S}$ is a singular point of its Hilbert scheme in most cases, Theorem ~\ref{B} (1) also contrasts several earlier results on ribbons and carpets on curves and surfaces respectively:  B\u{a}nic\u{a} and Manolache, in \cite{BM}, proved that the Hilbert points of ribbons in $\mathbb{P}^3$ supported on conics are smooth; Bayer and Eisenbud, in \cite{BE}, proved that the Hilbert
points of canonically embedded ribbons on $\mathbb{P}^1$ are smooth; and Gonz\'alez proved in \cite{Gon}, the smoothness of the Hilbert point for most ribbons on curves of arbitrary genus. It was also shown in \cite{GGP} and \cite{MR} that K3 carpets on Enriques surfaces and regular K-trivial ribbons on Enriques manifolds represent smooth points of the corresponding Hilbert schemes.
  
\subsection*{Acknowledgements} The first author thanks David Eisenbud for a conversation in Mysore, India in December of 2019, during which it became clear that one should prove much more general results than those in \cite{GP97} and that it would be interesting to do so. The conversation and results in \cite{BGG} motivated this article. We would also like to thank Jean-Marc Dr\'ezet for pointing us out an error in the statement of  Theorem ~\ref{countingp2} where the unique abstract non-split K3 carpet on $\mathbb{P}^2$ had been stated to be projective whereas it should have been non-projective as our proof shows. We are grateful to the anonymous referee for many useful suggestions and corrections that substantially helped in the improvement of this exposition.

\subsection*{Convention} We will always work over the complex numbers $\mathbb{C}$. For a smooth variety $X$, $K_X$ denotes the canonical bundle of $X$. The symbol ``$\sim$'' stands for linear equivalence and the symbol ``$\equiv$'' stands for numerical equivalence of line bundles or divisors. 


\section{Abstract and embedded K3 carpets}
This section is to investigate the existence of K3 carpets supported on a minimal rational surface $S$. We will deal with two cases, namely when $S=\mathbb{F}_e$ and when $S=\mathbb{P}^2$. We start with some basic things on ribbons and carpets. 
\subsection{Ribbons and K3 carpets} We give the precise definition of ribbon below.
\begin{definition}\phantomsection\label{defropes}
Let $Y$ be a reduced connected scheme and let $\mathcal{E}$ be a line bundle on Y. A ribbon on Y with conormal bundle $\mathcal{E}$ is a scheme $\widetilde{Y}$ with $\widetilde{Y}_{red}=Y$ such that $\mathcal{I}_{Y/Y'}^2=0$, and $\mathcal{I}_{Y/Y'}=\mathcal{E}$ as $\mathcal{O}_Y$ modules.
\end{definition}
The following properties of ribbons (more generally for ropes) were proven in \cite{BE} and \cite{Gon}.
\begin{theorem}\phantomsection\label{charropes}
Let $Y$ be a reduced connected scheme and $\mathcal{E}$ be a line bundle on Y.
\begin{itemize}[leftmargin=0.3in]
    \item[$(1)$] A ribbon $\widetilde{Y}$ with conormal bundle $\mathcal{E}$ is defined by an element $[e_{\widetilde{Y}}]\in \textrm{Ext}^1(\Omega_Y,\mathcal{E})$. The ribbon is split if and only if $[e_{\widetilde{Y}}]=0$.
\end{itemize}
Assume further that $Y$ is a smooth variety and $i:Y\hookrightarrow Z$ is a closed immersion into another smooth variety $Z$.
\begin{itemize}[leftmargin=0.3in]
    \item[$(2)$] There is an one-to-one correspondence between pairs $(\widetilde{Y},\widetilde{i})$ where $\widetilde{Y}$ is a ribbon with conormal bundle $\mathcal{E}$ and $\widetilde{i}:\widetilde{Y}\to Z$ is a morphism extending $i:Y\hookrightarrow Z$ and elements $\tau\in \textrm{Hom}(\mathcal{N}_{Y/Z}^*,\mathcal{E})$.
    \item[$(3)$] If $\tau\in \textrm{Hom}(\mathcal{N}_{Y/Z}^*,\mathcal{E})$ corresponds to $(\widetilde{Y},\widetilde{i})$, the $\widetilde{i}$ is an embedding if and only if $\tau$ is surjective.
    \item[$(4)$] If $\tau\in \textrm{Hom}(\mathcal{N}_{Y/Z}^*,\mathcal{E})$ corresponds to $(\widetilde{Y},\widetilde{i})$ then $\tau$ is mapped by the connecting homomorphism onto $[e_{\widetilde{Y}}]$.
\end{itemize}
\end{theorem}
K3 carpets were defined as ribbons on surfaces satisfying some additional properties (see \cite{GGP}, Definition 1.2).
\begin{definition}\label{k3carpets}
A K3 carpet $\widetilde{S}$ on a smooth regular surface $S$ is a ribbon on $S$ such that the dualizing sheaf $K_{\widetilde{S}}$ is trivial and $H^1(\mathcal{O}_{\widetilde{S}})=0$.
\end{definition}
The following characterizing lemma was proven in \cite{GP97}, Proposition $1.5$.
\begin{lemma}\label{chark3}
Let $S$ be a smooth regular surface and let $\widetilde{S}$ be a ribbon supported on $S$. Let $\mathcal{L}$ be the dual of the ideal sheaf defining $S$ in $\widetilde{S}$. The $\widetilde{S}$ is a K3 carpet if and only if $\mathcal{L}\cong K_S^*$.
\end{lemma} 
We refer to \cite{MR} for a generalization of K3 carpets in higher dimensions.
\subsection{K3 carpets on Hirzebruch surfaces} Recall that the Hirzebruch surfaces $S=\mathbb{F}_e$ ($e\geq 0$) are, by definition, the projective bundles $\mathbb{P}(\mathcal{E})$, where $\mathcal{E}=\mathcal{O}_{\mathbb{P}^1}\oplus\mathcal{O}_{\mathbb{P}^1}(-e)$. Let $\pi':S\to\mathbb{P}^1$ be the natural morphism. The line bundles on $S$ are of the form $\mathcal{O}_S(aC_0+bf)$ where $\mathcal{O}_S(C_0)=\mathcal{O}_{\mathbb{P}(\mathcal{E})}(1)$ and $\mathcal{O}_S(f)=\pi'^*\mathcal{O}_{\mathbb{P}^1}(1)$. The line bundle $\mathcal{O}_S(aC_0+bf)$ is very ample if and only if $b\geq ae+1$.
\begin{lemma}\label{cohomology}
Suppose $j: S \hookrightarrow \mathbb{P}^n$ be an embedding of a Hirzebruch surface $S$ induced by the complete linear series of $\mathcal{O}_S(1) = \mathcal{O}_S(aC_0 + bf)$ (and hence $b \geq ae + 1$ and $n+1 = h^0(\mathcal{O}_S(1))$). Let $N+1 \geq h^0(\mathcal{O}_S(1))$ and let $k: \mathbb{P}^n \hookrightarrow \mathbb{P}^N$ be a linear embedding. Let $i = k \circ j$. Let $\pi': \mathbb{F}_e \to \mathbb{P}^1$ be the natural fibration. Then
\begin{itemize}
    \item[$(1)$] $h^0(\pi'^*T_{\mathbb{P}^1} \otimes K_S) = 0, h^1(\pi'^*T_{\mathbb{P}^1} \otimes K_S) = 1, h^2(\pi'^*T_{\mathbb{P}^1} \otimes K_S) = 0$. 
\item[$(2)$] $h^0(T_{S/\mathbb{P}^1} \otimes K_S) = 0, h^1(T_{S/\mathbb{P}^1} \otimes K_S) = 1, h^2(T_{S/\mathbb{P}^1} \otimes K_S) = 0$.
\item[$(3)$] $h^0(T_{\mathbb{P}^N}|_S \otimes K_S) = (N+1)h^0((aC_0 + bf) \otimes K_S), h^1(T_{\mathbb{P}^N}|_S \otimes K_S) = 1,  h^2(T_{\mathbb{P}^N}|_S \otimes K_S) = 0$.
\item[$(4)$] $h^0(T_S \otimes K_S) = h^2(T_S \otimes K_S) = 0, h^1(T_S \otimes K_S) = 2 $.
\end{itemize}
\end{lemma}

\noindent\textit{Proof.} Recall that $K_S = \mathcal{O}_S(-2C_0 -(e+2)f)$ is the canonical bundle of $S$.\par 
($1$) Note that $\pi'^*T_{\mathbb{P}^1} = \mathcal{O}_S(2f)$, thus $\pi'^*T_{\mathbb{P}^1} \otimes K_S = \mathcal{O}_S(2f) \otimes K_S$. Now, $h^i( \mathcal{O}_S(2f) \otimes K_S) = h^{2-i}(\mathcal{O}_S(-2f))$ by Serre duality. The assertions follow by using $\pi'_*(\mathcal{O}_S) = \mathcal{O}_{\mathbb{P}^1}$ and the projection formula.\par 

($2$) We have $T_{S/\mathbb{P}^1} = \mathcal{O}_S(2C_0 + ef)$. Hence $T_{S/\mathbb{P}^1} \otimes K_S = \mathcal{O}_S(-2f)$. The assertions again follow by the projection formula and by the fact that $\pi'_*(\mathcal{O}_S) = \mathcal{O}_{\mathbb{P}^1}$.\par 
($3$) We tensor the restriction of the Euler sequence to $S$ by $K_S$ to get the following
\begin{equation*}
    0 \to K_S \to \mathcal{O}_S(aC_0+bf)^{\oplus N+1} \otimes K_S \to T_{\mathbb{P}^N}|_S \otimes K_S \to 0.
\end{equation*}
Since $\mathcal{O}_S(aC_0 + bf)$ is very ample, we have that $h^1((aC_0 + bf) \otimes K_S) = h^2((aC_0 + bf) \otimes K_S) = 0$ by Kodaira vanishing. One checks that $h^0(K_S) = 0$, $h^1(K_S) = 0$ and  $h^2(K_S)= 1$. The assertions follow by taking the long exact sequence associated to the above short exact sequence.\par 

($4$) The assertions follow from the long exact sequence associated to
\begin{equation*}
    0 \to T_{S/\mathbb{P}^1} \otimes K_S \to T_S \otimes K_S \to \pi'^*(T_{\mathbb{P}^1}) \otimes K_S \to 0.
\end{equation*}
and the previous parts.\QEDB\par 

\vspace{5pt}

We are interested in counting the number of abstract K3 carpets on Hirzebruch surfaces, as well as the number of embedded K3 carpets on the embedded Hirzebruch surfaces. It follows from Theorem ~\ref{charropes} and Lemma ~\ref{chark3} that we need the dimension of the cohomology group $H^1(T_S\otimes K_S)$ for the first one, and the dimension of $H^0(\mathcal{N}_{S/\mathbb{P}^N}\otimes K_S)$ for the second one. To compute their values, we need the following lemma. 
\begin{lemma}\label{higher direct images} In the situation of Lemma ~\ref{cohomology}, we have the following
 \begin{itemize}
     \item[$(1)$] $\pi'_*(T_{S/\mathbb{P}^1} \otimes K_S) = \mathcal{O}_{\mathbb{P}^1}(-2)$,  $R^1\pi'_*(T_{S/\mathbb{P}^1} \otimes K_S) = 0$.
     \item[$(2)$] $\pi'_*(\pi'^*(T_{\mathbb{P}^1}) \otimes K_S) = 0$,  $R^1\pi'_*(\pi'^*(T_{\mathbb{P}^1}) \otimes K_S) = \mathcal{O}_{\mathbb{P}^1}$.
     \item[$(3)$] $\pi'_*(T_S \otimes K_S) = \mathcal{O}_{\mathbb{P}^1}(-2)$,  $R^1\pi'_*(T_S \otimes K_S) = \mathcal{O}_{\mathbb{P}^1} $.
     \item[$(4)$] There is an exact sequence 
     \begin{equation*}
    0 \to \pi'_*(\mathcal{O}_S((aC_0 + bf) \otimes K_S)^{\oplus N+1})  \to \pi'_*(T_{\mathbb{P}^N}|_S \otimes K_S) \to R^1\pi'_*(K_S) = \mathcal{O}_{\mathbb{P}^1}(-2)  \to 0,
     \end{equation*}
     and $R^1\pi'_*(T_{\mathbb{P}^N}|_S \otimes K_S) = 0$.
     \item[$(5)$] There is an exact sequence 
     \begin{equation*}
    0 \to \mathcal{O}_{\mathbb{P}^1}(-2) \to \pi'_*(T_{\mathbb{P}^N}|_S \otimes K_S) \to \pi'_*(\mathcal{N}_{S/\mathbb{P}^N} \otimes K_S) \to \mathcal{O}_{\mathbb{P}^1} \to 0,
     \end{equation*}
     and $R^1\pi'_*(\mathcal{N}_{S/\mathbb{P}^N} \otimes K_S) = 0$.
      
 \end{itemize}
\end{lemma}
\noindent\textit{Proof.} ($1$) and ($2$) are clear since $T_{S/\mathbb{P}^1} \otimes K_S \cong \mathcal{O}_S (-2f)$ and $\pi'^*T_{\mathbb{P}^1} \otimes K_S = \mathcal{O}_S (-2C_0 -ef)$. 
($3$) follows using ($1$) and ($2$) and applying $\pi'_*$ to the exact sequence 
\begin{equation*}
    0 \to \mathcal{O}_S (-2f) \to T_S \otimes K_S \to \mathcal{O}_S (-2C_0 -ef) \to 0.
\end{equation*} 
($4$) follows by applying $\pi'_*$ to the sequence obtained by tensoring the Euler sequence by $K_S$ to obtain
\begin{equation*}
    0 \to K_S \to \mathcal{O}_S((aC_0 + bf) \otimes K_S)^{\oplus N+1}) \to T_{\mathbb{P}^N}|_S \otimes K_S \to 0
\end{equation*} 
and relative Serre duality.
($5$) follows by applying $\pi'_*$ to the following exact sequence 
\begin{equation*}
  0 \to T_S \otimes K_S \to T_{\mathbb{P}^N}|_S \otimes K_S \to \mathcal{N}_{S/\mathbb{P}^N} \otimes K_S \to 0 
\end{equation*}
and the previous parts.\QEDB\par 

\vspace{5pt}

We have seen in the previous lemma that $R^1\pi'_*(\mathcal{N}_{S/\mathbb{P}^N} \otimes K_S) = 0$. Thus, the dimensions of cohomology groups of $\mathcal{N}_{S/\mathbb{P}^N} \otimes K_S$ can be calculated by pushing it forward to $\mathbb{P}^1$. The following proposition gives an exact sequence where $\pi'_*(\mathcal{N}_{S/\mathbb{P}^N} \otimes K_S)$ fits, and we use that exact sequence to calculate $h^0(\pi'_*(\mathcal{N}_{S/\mathbb{P}^N} \otimes K_S))$.

\begin{lemma}\label{h1} In the situation of Lemma ~\ref{cohomology}, we have the following;
\begin{itemize}
     \item[$(1)$] There is an exact sequence
     \begin{equation*}
    0 \to \pi'_*(\mathcal{O}_S((aC_0 + bf) \otimes K_S)^{\oplus N+1})  \to \pi'_*(\mathcal{N}_{S/\mathbb{P}^N} \otimes K_S) \to \mathcal{O}_{\mathbb{P}^1} \to 0.
     \end{equation*}
     \item[$(2)$] $h^0(\mathcal{N}_{S/\mathbb{P}^N} \otimes K_S) = (N+1)\frac{(a-1)}{2}(2b-2-ae)$. In particular if the embedding $i$ is induced by the complete linear series of $\mathcal{O}_S(aC_0+bf)$, then $N+1 = h^0(\mathcal{O}_S(aC_0+bf)$ and $h^0(\mathcal{N}_{S/\mathbb{P}^N} \otimes K_S) = \frac{1}{4}(a^2-1)((2b-ae)^2 - 4) + 1.$ We also have $h^1(\mathcal{N}_{S/\mathbb{P}^N} \otimes K_S) = h^2(\mathcal{N}_{S/\mathbb{P}^N} \otimes K_S) = 0$.
\end{itemize}
\end{lemma}

\noindent\textit{Proof.} ($1$)  
Choose an injective morphism $\psi: \mathcal{O}_S(C_0+bf) \to \mathcal{O}_S(aC_0+bf)$. Let $h^0(\mathcal{O}_S(C_0+bf)) = m+1$. This induces the commutative diagram below.
\[ 
\begin{tikzcd}
     0 \arrow{r} & K_S \arrow{r} \arrow[d, "\cong"] &(\mathcal{O}_S(C_0 + bf) \otimes K_S ) \otimes H^0(\mathcal{O}_S(C_0+bf) ) \arrow{d} \arrow{r} & T_{\mathbb{P}^m}|_S \otimes K_S \arrow{r} \arrow[d, "\alpha"] & 0 \\
     0 \arrow{r} & K_S  \arrow{r} &(\mathcal{O}_S(aC_0 + bf) \otimes K_S )^{\oplus N+1}  \arrow{r} & T_{\mathbb{P}^N}|_S \otimes K_S \arrow{r} & 0 
\end{tikzcd} \]
Taking $\pi'_*$ and noting that $R^i\pi'_*((\mathcal{O}_S(C_0 + bf) \otimes K_S ) \otimes H^0(\mathcal{O}_S(C_0+bf) )) = 0$ for $i = 1$, $2$ and using Lemma ~\ref{higher direct images} (4), we have the following commutative diagram. 
\[ 
\begin{tikzcd}
    0 \arrow{r} & 0 \arrow{r} \arrow{d} & \pi'_*(T_{\mathbb{P}^m}|_S \otimes K_S) \arrow[d, "\phi=\pi'_*\alpha"] \arrow[r, "f"] &  R^1\pi'_*K_S \arrow[d, "\cong"] \arrow{r} & 0  \\
     0 \arrow{r} & \pi'_*((\mathcal{O}_S(aC_0 + bf) \otimes K_S )^{\oplus N+1}  \arrow{r}  & \pi'_*(T_{\mathbb{P}^N}|_S \otimes K_S) \arrow[r,"h"]  &  R^1\pi'_*K_S \arrow{r} & 0 
\end{tikzcd} \]
\textcolor{black}{The right commutative square gives a section of $h$ and hence the bottom exact sequence splits and} we have that $$\phi: \pi'_*(T_{\mathbb{P}^m|_S} \otimes K_S) \to \pi'_*(T_{\mathbb{P}^N}|_S \otimes K_S) = R^1\pi'_*K_S \oplus \pi'_*((\mathcal{O}_S(aC_0 + bf) \otimes K_S )^{\oplus N+1} $$ and $\phi$ induces an isomorphism between $\pi'_*(T_{\mathbb{P}^m|_S} \otimes K_S)$ and $R^1\pi'_*K_S$.\par 

Now consider the commutative diagram induced by $\alpha$.
\[ 
\begin{tikzcd}
     0 \arrow{r} & T_S \otimes K_S \arrow{r} \arrow[d, "\cong"] &T_{\mathbb{P}^m}|_S \otimes K_S \arrow[d, "\alpha"] \arrow{r} & \mathcal{N}_{S/\mathbb{P}^m} \otimes K_S \arrow{r} \arrow{d} & 0 \\
     0 \arrow{r} & T_S \otimes K_S  \arrow{r} & T_{\mathbb{P}^N}|_S \otimes K_S  \arrow{r} & \mathcal{N}_{S/\mathbb{P}^N} \otimes K_S \arrow{r} & 0 
\end{tikzcd} \]
Applying $\pi'_*$ and using Lemma ~\ref{higher direct images} we get 
\[ 
\begin{tikzcd}
     0 \arrow{r} & \pi'_*(T_S \otimes K_S) \arrow{r} \arrow[d, "\cong"] &\pi'_*(T_{\mathbb{P}^m|_S} \otimes K_S) \arrow[d, "\phi"] \arrow{r} & \pi'_*(\mathcal{N}_{Y/\mathbb{P}^m} \otimes K_S) \arrow{r} \arrow{d} & R^1\pi'_*(T_S \otimes K_S) \arrow{r} \arrow{d} & 0 \\
     0 \arrow{r} & \pi'_*(T_S \otimes K_S)  \arrow[r, "p"] & \pi'_*(T_{\mathbb{P}^N}|_S \otimes K_S)  \arrow{r} & \pi'_*(\mathcal{N}_{S/\mathbb{P}^N} \otimes K_S) \arrow{r} & R^1\pi'_*(T_S \otimes K_S)  \arrow{r} & 0 
\end{tikzcd} \]
 Since the leftmost vertical map is an isomorphism we have that $p$ maps $\pi'_*(T_S \otimes K_S)$ isomorphically onto $\phi(\pi'_*(T_{\mathbb{P}^m}|_S \otimes K_S)) = R^1\pi'_*K_S$. By the previous paragraph and Lemma ~\ref{higher direct images} (3), we have the exact sequence 
 \begin{equation*}
    0 \to \pi'_*(\mathcal{O}_S((aC_0 + bf) \otimes K_S)^{\oplus N+1})  \to \pi'_*(\mathcal{N}_{S/\mathbb{P}^N} \otimes K_S) \to \mathcal{O}_{\mathbb{P}^1} \to 0.
     \end{equation*}
\indent $(2)$ Notice that we have $R^1\pi'_*(\mathcal{N}_{S/\mathbb{P}^N} \otimes K_S) = 0$ by Lemma ~\ref{higher direct images} (5) and hence $$h^i(\mathcal{N}_{S/\mathbb{P}^N} \otimes K_S) = h^i(\pi'_*(\mathcal{N}_{S/\mathbb{P}^N} \otimes K_S)).$$ 

Thus, it follows from the exact sequence of $(1)$ that, $$h^0(\mathcal{N}_{S/\mathbb{P}^N} \otimes K_S) = (N+1)h^0(\pi'_*(\mathcal{O}_S((aC_0 + bf) \otimes K_S))) + 1$$ and $h^1(\mathcal{N}_{S/\mathbb{P}^N} \otimes K_S) = h^2(\mathcal{N}_{S/\mathbb{P}^N} \otimes K_S)=0$ (to see the vanishings, use  
$R^1\pi'_*(\mathcal{O}_S((aC_0 + bf) \otimes K_S)$=0,  
and Kodaira vanishing).
 
The assertion follows from the fact that 
\begin{equation*}
    h^0(\mathcal{O}_S((aC_0+bf) \otimes K_S)) = \frac{(a-1)}{2}(2b-2-ae)\textrm{ and } h^0(\mathcal{O}_S(aC_0+bf)) = \frac{(a+1)}{2}(2b+2-ae).
\end{equation*}
That finishes the proof of the proposition.\QEDB\par 

\section{Projective and non-projective K3 carpets}

In contrast to ribbons on curves, not all carpets are projective, even if all of them are proper
or even if , as is the case with rational surfaces, they are supported on a projective surface. A natural question to ask 
about K3 carpets on minimal rational surfaces is whether there exist families of projective K3 carpets. We show this is true 
in Theorem ~\ref{counting carpets} below. Next step is to describe the loci parametrizing K3 carpets on a given minimal rational surface, which 
we do as well. 
 
As we will see, the situation resembles that of smooth K3 surfaces. 
We also explore embedded K3 carpets, which naturally are dependent on the embeddings of rational surfaces in projective
space. We give precise results on all of the above in the following theorem:

\begin{theorem}\label{counting carpets} In the situation of Lemma ~\ref{cohomology}, we have the following;
  \begin{itemize}
      \item[$(1)$] The abstract non-split $K3$ carpets on $S$ are parametrized by the projective space of lines in the vector space $H^1(T_S \otimes K_S)$ of dimension $h^1(T_S \otimes K_S) = 2$. Hence the non-split abstract $K3$ carpets on $S$ are parametrized by the projective line $\mathbb{P}^1$.
     \item[$(2)$] The non-split projective $K3$ carpets on $S$ are parametrized by a non-empty countable subset of the projective line . They are in $1-1$ correspondence to the classes of primitive ample divisors in $NS(S)$.
     \item[$(3)$] The $K3$ carpets on $S$ with an embedding inside $\mathbb{P}^N$ that extends the embedding $i$ are parametrized by a non-empty Zariski open subset of the  projective space of lines in the vector space $H^0(\mathcal{N}_{S/\mathbb{P}^N} \otimes K_S)$. The Zariski open set is the projectivization of of the open set of $H^0(\mathcal{N}_{S/\mathbb{P}^N} \otimes K_S)$ consisting of nowhere vanishing sections. Moreover, $$h^0(\mathcal{N}_{S/\mathbb{P}^N} \otimes K_S) = (N+1)\frac{(a-1)}{2}(2b-2-ae) + 1.$$ 
  \end{itemize}

\end{theorem}
\noindent\textit{Proof.} (1) The assertion is clear, $h^1(T_S\otimes K_S)$ has been calculated in Lemma ~\ref{cohomology}.\par 
(2) Follow the proof of \cite{GGP}, Theorem 2.5 word by word.


(3)  Using \cite{Gon}, in order to prove the statement,
we need to show that the Zariski open subset of the  projective space of lines in the vector space $H^0(\mathcal{N}_{S/\mathbb{P}^N} \otimes K_S)$ consisting of nowhere vanishing sections is non-empty. In other words we want to show that there exist at least one nowhere vanishing section of the vector bundle $\mathcal{N}_{S/\mathbb{P}^N} \otimes K_S$. Consider the exact sequence  
\begin{equation*}
    0 \to K_S  \to (\mathcal{O}_S((aC_0 + bf) \otimes K_S)^{\oplus N+1}) \to T_{\mathbb{P}^N}|_S \otimes K_S \to 0.
     \end{equation*}
 \indent Notice that $\mathcal{O}_S(aC_0+bf) \otimes K_S = \mathcal{O}_S((a-2)C_0+(b-e-2)f)$ is globally generated. Indeed, since $b \geq ae+1$, we have $$b-e-2 \geq ae+1-e-2 = ae-e-1 \geq ae-2e$$ for $e \geq 1$ and hence $\mathcal{O}_S(aC_0+bf) \otimes K_S$ is base point free for $e\geq 1$. For $e =0$, $S = \mathbb{P}^1 \times \mathbb{P}^1$, hence if $a \geq 2$ and $b \geq 2$ we still have $\mathcal{O}_S(aC_0+bf) \otimes K_S$ is base point free.\par 
 
 Assume $a \geq 2$. From Lemma ~\ref{cohomology}, we have that $$h^0(T_{\mathbb{P}^N}|_S \otimes K_S) = (N+1)h^0(\mathcal{O}_S((aC_0 + bf) \otimes K_S)) = (N+1) \frac{(a-1)}{2}(2b-2-ae)$$ which is non-zero since $a \geq 2$. Thus, $T_{\mathbb{P}^N}|_S \otimes K_S$ is globally generated, since $\mathcal{O}_S((aC_0 + bf) \otimes K_S)^{\oplus N+1}$ is base point free and it surjects onto $T_{\mathbb{P}^N}|_S \otimes K_S$. Consider the sequence    
 \begin{equation*}
    0 \to T_S \otimes K_S  \to  T_{\mathbb{P}^N}|_S \otimes K_S \to \mathcal{N}_{S/\mathbb{P}^N} \otimes K_S \to 0.
     \end{equation*}
We already have that $h^0(\mathcal{N}_{S/\mathbb{P}^N} \otimes K_S) = (N+1)\frac{(a-1)}{2}(2b-2-ae)+1$ and hence $h^0(\mathcal{N}_{S/\mathbb{P}^N} \otimes K_S) \neq 0$.  
Thus, we have that $\mathcal{N}_{S/\mathbb{P}^N} \otimes K_S$ is base point free. Now,  rank$(\mathcal{N}_{S/\mathbb{P}^N} \otimes K_S) = N-2$. It is easy to see using $N+1 \geq \frac{1}{2}(a+1)(2b+2-ae)$ that, 
$N-2 > 2$ and hence $H^0(\mathcal{N}_{S/\mathbb{P}^N} \otimes K_S)$ contains a nowhere vanishing section. \par
It remains to show the existence of a nowhere vanishing section when $a = 1$.  
But this follows from \cite{GP97}, Proposition 1.7.  
That ends the proof. \QEDB

\vspace{5pt}

We thank the referee for suggesting the following remark that generalizes \cite{GGP}, Theorem 2.5 and Theorem ~\ref{counting carpets} (2).
\begin{remark}
 Let $S$ be a smooth surface and assume that $S$ supports more than one (abstract) non-split K3 carpet. Following the proof of \cite{GGP2}, Theorem 2.5, one can prove in this case, that the projective non-split K3 carpets on $S$ are parametrized by a countable union of hyperplanes in $\mathbb{P}(H^1(T_S\otimes K_S))$.
\end{remark}

\subsection{K3 carpets on $\mathbb{P}^2$} Now we count the number of abstract K3 carpets on $\mathbb{P}^2$ and the number of embedded K3 carpets on the $d$-uple embeddings of $\mathbb{P}^2$.
\begin{lemma}\label{p2} Suppose that $j: S = \mathbb{P}^2 \hookrightarrow \mathbb{P}^n$ be an embedding induced by the complete linear series $\mathcal{O}_{\mathbb{P}^2}(d)$ (and hence $n+1 = h^0(\mathcal{O}_{\mathbb{P}^2}(d))$). Let $N+1 \geq h^0(\mathcal{O}_{\mathbb{P}^2}(d))$ and let $k: \mathbb{P}^n \hookrightarrow \mathbb{P}^N$ be a linear embedding. Let $i = k \circ j$. 
\begin{itemize}
\item[$(1)$] 
     Then there is an exact sequence,
     \begin{equation*}
         0 \to \mathcal{O}_{\mathbb{P}^2}^{\oplus 3}(1) \to \mathcal{O}_{\mathbb{P}^2}^{\oplus N+1}(d) \to \mathcal{N}_{S/\mathbb{P}^N} \to 0.
     \end{equation*}
     \item[$(2)$] $h^0(\mathcal{N}_{S/\mathbb{P}^N} \otimes K_S) = \frac{1}{2}(N+1)(d-1)(d-2)$. In particular if the embedding $i$ is induced by the complete linear series of $\mathcal{O}_{\mathbb{P}^2}(d)$, then  $ h^0(\mathcal{N}_{S/\mathbb{P}^N} \otimes K_S) = \frac{1}{4}(d+2)(d+1)(d-1)(d-2)$. We also have $h^1(\mathcal{N}_{S/\mathbb{P}^N} \otimes K_S) = h^2(\mathcal{N}_{S/\mathbb{P}^N} \otimes K_S) = 0$
     \end{itemize}
\end{lemma}
\noindent\textit{Proof.} (1) The following commutative diagram with exact rows and columns shows $\mathcal{F}\cong\mathcal{N}_{S/\mathbb{P}^N}$. 

\adjustbox{scale=0.8,center}{
\begin{tikzcd}
 & 0 \arrow{d} & 0 \arrow{d} & & \\
  0 \arrow{r} &  \mathcal{O}_{\mathbb{P}^2} \arrow{r} \arrow{d} &  \mathcal{O}_{\mathbb{P}^2} \arrow{r} \arrow{d} & 0 \arrow{d}  \arrow{r} & 0\\
 0 \arrow{r} & \mathcal{O}_{\mathbb{P}^2}^{\oplus 3}(1) \arrow{r} \arrow{d} & \mathcal{O}_{\mathbb{P}^2}^{\oplus N+1}(d) \arrow{r} \arrow{d} & \mathcal{F} \arrow{d} \arrow{r} & 0  \\
   0 \arrow{r} & T_{\mathbb{P}^2} \arrow{r} \arrow{d} & T_{\mathbb{P}^N}|_S \arrow{r} \arrow{d} & \mathcal{N}_{S/\mathbb{P}^N}\arrow{d} \arrow{r} & 0\\ 
  & 0 & 0 & 0 &  
\end{tikzcd} 
}

\indent (2) The assertion follows by tensoring the exact sequence of the previous part with $K_S$ and using cohomology calculations of the long exact sequence associated to it. \QEDB
 \begin{theorem}\label{countingp2} In the situation of Lemma ~\ref{p2}, we have the following.
 \begin{itemize}

    \item[$(1)$] There exists a unique abstract non-split $K3$ carpet on $S$ which is non-projective. 
    \item[$(2)$] The $K3$ carpets on $S$ with an embedding inside  $\mathbb{P}^N$ that extends the embedding $i$ are parametrized by a non-empty Zariski open subset of the  projective space of lines in the vector space $H^0(\mathcal{N}_{S/\mathbb{P}^N} \otimes K_S)$. The Zariski open set is the projectivization of of the open set of $H^0(\mathcal{N}_{S/\mathbb{P}^N} \otimes K_S)$ consisting of nowhere vanishing sections. Moreover, $$h^0(\mathcal{N}_{S/\mathbb{P}^N} \otimes K_S) = \frac{1}{2}(N+1)(d-1)(d-2).$$ In particular, there do not exist an embedded $K3$ carpet for the embedding given by $d =2$, i.e. for the second Veronese embedding.
    \end{itemize}
\end{theorem}
\noindent\textit{Proof.} The proof follows exactly along the same lines of the proof of Theorem ~\ref{counting carpets}.  (1) follows since $h^1(T_S\otimes K_S)=1$ and by Poincare duality $\int_D$ is a non-degenerate form on $H^1(T_S \otimes K_S)$ 
(alternatively, it follows from \cite{Dr}, Theorem $8.0.2$). For (2), notice that  $\mathcal{O}_{\mathbb{P}^2}^{\oplus N+1}(d-3)$ is base point free for $d \geq 3$ and there is an exact sequence 
  \begin{equation}\label{3.1}
     0 \to \mathcal{O}_{\mathbb{P}^2}^{\oplus 3}(-2)  \to \mathcal{O}_{\mathbb{P}^2}^{\oplus N+1}(d-3) \to \mathcal{N}_{S/\mathbb{P}^N} \otimes K_S \to 0.
 \end{equation}
obtained by tensoring the exact sequence of Lemma ~\ref{p2}, (1) by $K_S=\mathcal{O}_{\mathbb{P}^2}(-3)$. which implies $\mathcal{N}_{S/\mathbb{P}^N} \otimes K_S$ is base point free Since rank$(\mathcal{N}_{S/\mathbb{P}^N} \otimes K_S) > 2$,  a general section of $\mathcal{N}_{S/\mathbb{P}^N} \otimes K_S$ is nowhere vanishing. \QEDB\par

\section{Smoothings of abstract and embedded K3 carpets}
We study the smoothings of the K3 carpets constructed in the previous section. We start with the following result (see \cite{Gon}, Lemma 3.3 and Proposition 3.7) that is crucial for us.
\begin{lemma}\label{deformation}
Let $X$ be a smooth projective variety and $\varphi: X \to \mathbb{P}^N$ be a morphism. Suppose that $\varphi = i \circ \pi$ where $\pi: X \to Y$ be a finite flat morphism with $Y$ smooth and trace zero module $\mathcal{E}$ and $i: Y \hookrightarrow \mathbb{P}^N$ be an embedding. Let $\mathcal{N}_{\pi}, \mathcal{N}_{\varphi}, \mathcal{N}_{Y/\mathbb{P}^N}$ be the normal sheaves of the morphisms $\varphi, \pi,$ and $i$ respectively. Then 

    there exists an exact sequence 
    \begin{equation*}
    0 \to \mathcal{N}_{\pi}  \to  \mathcal{N}_{\varphi} \to \pi^*\mathcal{N}_{Y/\mathbb{P}^N} \to 0.
     \end{equation*}

Taking cohomology and using projection formula we have the following exact sequence 
\begin{equation*}
    0 \to H^0(\mathcal{N}_{\pi}) \to H^0(\mathcal{N}_{\varphi})  \xrightarrow{\nu_1\oplus\nu_2} H^0(\mathcal{N}_{Y/\mathbb{P}^N}) \oplus H^0(\mathcal{N}_{Y/\mathbb{P}^N} \otimes \mathcal{E}) \to H^1(\mathcal{N}_{\pi})
     \end{equation*}
 where $\nu_1:H^0(\mathcal{N}_{\varphi}) \to H^0(\mathcal{N}_{Y/\mathbb{P}^N}$ and $\nu_2 : H^0(\mathcal{N}_{\varphi}) \to H^0(\mathcal{N}_{Y/\mathbb{P}^N}) \otimes \mathcal{E})$ are the induced maps.

\end{lemma}

The following proposition is fundamental in the study of the deformation theory of morphisms and was developed by Gallego, Gonz\'alez and the first author. We will use this result throughout this section.
\begin{proposition}\label{main}
In the situation of Lemma ~\ref{deformation}, let $\Phi:\mathcal{X}\to\mathbb{P}^N_T$ be a deformation of $\varphi$ along a smooth pointed affine curve $(T,0)$. Assume that 
\begin{itemize}
    \item[$(a)$] $\mathcal{X}$ is irreducible and reduced;
    \item[$(b)$] $\mathcal{X}_t$ is smooth, irreducible and projective for all $t\in T$;
    \item[$(c)$] $\mathcal{X}_0=X$ and $\Phi_0=\varphi$.
\end{itemize} 
Let $\Delta$ be the first infinitesimal neighborhood of $t\in T$. Let $\widetilde{X}=\mathcal{X}_{\Delta}$ and $\widetilde{\varphi}=\Phi_{\Delta}$. Then $\widetilde{\varphi}$ naturally defines an element (which we again call $\widetilde{\varphi}$) in $H^0(\mathcal{N}_{\varphi})$. If $\nu_2(\widetilde{\varphi}) \in H^0(\mathcal{N}_{Y/\mathbb{P}^N} \otimes \mathcal{E})$ is a nowhere vanishing section of the the vector bundle $\mathcal{N}_{Y/\mathbb{P}^N} \otimes \mathcal{E}$, then the following happens.
\begin{itemize}
    \item[$(1)$] After possibly shrinking $T$, we have that for $t \in T, t \neq 0$, the morphism $\Phi_t: \mathcal{X}_t \to \mathbb{P}^N_t $ is an embedding. 
    \item[$(2)$] The central fibre of the image, $($Im $\Phi)_0$ is a multiple structure on Im$(\varphi)$ with conormal bundle $\mathcal{E}$. Since Im$(\Phi)_t \cong \mathcal{X}_t$ and $\mathcal{X}_t$ is smooth, we have that Im$(\Phi) \to T$ is an embedded smoothing family of $($Im $\Phi)_0$.
\end{itemize}
\end{proposition}
\noindent\textit{Proof.} See \cite{GGP1}, Proposition $1.4$. The last assertion follows from equation $1.9$ in the proof of Proposition $1.4$, \cite{GGP1}. \QEDB\par 
\subsection{Smoothings of K3 carpets on Hirzebruch surfaces}  We now prove that projective K3 carpets are smoothable. We first deal with carpets on Hirzebruch surfaces $\mathbb{F}_e$ with $0\leq e\leq 2$. The reason we treat them first is because these surfaces admit smooth K3 double covers. We will use this result to show that the projective K3 carpets on the remaining Hirzebruch surfaces are smoothable as well, using a degeneration argument. Finally, we show smoothing results on projective $K3$ carpets on $\mathbb{P}^2$ embedded by an arbitrary very ample line bundle. 

\begin{theorem}\label{smoothing1}
 Let $\widetilde{S}$ be a projective $K3$ carpet supported on $S = \mathbb{F}_e$ with $0 \leq e \leq 2$. Then $\widetilde{S}$ is smoothable.
\end{theorem}

\noindent\textit{Proof.} Let $\widetilde{J}$ be a very ample line bundle on $\widetilde{S}$. Consider the embedding $\widetilde{i}: \widetilde{S} \hookrightarrow \mathbb{P}^N$ given by the complete linear series of $\widetilde{J}$. This induces an embedding $i: S \hookrightarrow \mathbb{P}^N$. Let $\widetilde{J}|_S = J$. Consider the following exact sequence: 
$$ 0 \to K_S \to \mathcal{O}_{\widetilde{S}} \to \mathcal{O}_S \to 0. $$
Tensoring by $\widetilde{J}$ we get 
$$ 0 \to K_S \otimes J \to \widetilde{J} \to J \to 0. $$

\textcolor{black}{Let $n = h^0(J)-1$. Since $H^1(K_S \otimes J) = 0$ we have that after a projection to a linear subspace $\mathbb{P}^n$ of codimension $h^0(K_S \otimes J)$, the composition of the embedding followed by the projection is induced by the complete linear series of $J := \mathcal{O}_S(aC_0 + bf)$ and hence satisfies the situation of Lemma ~\ref{cohomology}. Let $k: \mathbb{P}^n \hookrightarrow \mathbb{P}^N$ be the linear embedding.}

Now notice that the complete linear series of $|-2K_S|$ is base point free for $S=\mathbb{F}_e$ with $0 \leq e \leq 2$. By Bertini, we can choose a smooth curve section $C \in |-2K_S|$ and consider the double cover $\pi: X \to S$ branched along the smooth curve $C$. It follows that $K_X = \mathcal{O}_X$ and $H^1(\mathcal{O}_X) = H^1(\mathcal{O}_S) \oplus H^1(K_S) = 0$ since the trace zero module of $\pi$ is $K_S$. Thus, $X$ is a $K3$ surface. Now consider the composed morphism $\varphi: X \xrightarrow{\pi} S \xhookrightarrow{i} \mathbb{P}^N $. According to Lemma ~\ref{deformation}, we have an exact sequence 
\begin{equation*}
    0 \to H^0(\mathcal{N}_{\pi})  \to  H^0(\mathcal{N}_{\varphi}) \to H^0(\mathcal{N}_{Y/\mathbb{P}^N}) \oplus H^0(\mathcal{N}_{Y/\mathbb{P}^N} \otimes K_S) \to H^1(\mathcal{N}_{\pi})
     \end{equation*}
Observe that $\mathcal{N}_{\pi} = -2K_S|_C$ and it fits into
the following exact sequence 
\begin{equation}
0 \to \mathcal{O}_S \to -2K_S \to  -2K_S|_C \to 0.   
\end{equation}
Note that, since $H^1(-2K_S) = H^2(\mathcal{O}_S) = 0$ we have that $H^1(\mathcal{N}_{\pi}) = H^1(-2K_S|_C) = 0$. \par 
In Theorem ~\ref{counting carpets}, we showed that nowhere vanishing sections form a non-empty Zariski open subset in the projective space of lines of the vector space $\mathcal{N}_{S/\mathbb{P}^N} \otimes K_S$ and that $\widetilde{S}$ corresponds to such a nowhere vanishing section. This combined with $H^1(\mathcal{N}_{\pi}) = 0$ allows us to choose a first order deformation $\widetilde{\varphi} \in H^0(\mathcal{N}_{\varphi})$ i.e. a deformation $\widetilde{\varphi}: \mathcal{\widetilde{X}} \to \mathbb{P}_{\Delta}^N$ over the ring of dual numbers $\Delta$, which maps to $\widetilde{S} \in H^0(\mathcal{N}_{S/\mathbb{P}^N} \otimes K_S)$  under the map $\nu_2$ as defined in Lemma ~\ref{deformation}.\par 
Let $\widetilde{L} = \widetilde{\varphi}^*(\mathcal{O}_{\mathbb{P}_{\Delta}^N}(1))$ and $L = \varphi^*(\mathcal{O}_{\mathbb{P}^N}(1))$. Then the pair $(\widetilde{X}, \widetilde{L})$ is a first order deformation of the pair $(X,L)$. Since polarized $K3$ surfaces have a $19$ dimensional smooth algebraic formally universal deformation space $D$, one can choose a smooth affine algebraic curve $T$ in $D$ passing through the point $(X,L)$ with tangent vector given by $(\widetilde{X}, \widetilde{L})$. Pulling back the universal family on $D$ along $T$, we get a pair $(\mathcal{X} \xrightarrow{f} T, \mathcal{L})$ where $\mathcal{X}$ is smooth, irreducible and projective (since $L$ is ample) after possibly shrinking $T$ and $\mathcal{L}$ is a line bundle on $\mathcal{X}$. Since $T$ is affine \textcolor{black}{and smooth, we have that $f_*(\mathcal{L})$ is locally free of rank $h^0(L)$}, we have that $\mathcal{L}$ induces a morphism $\Psi: \mathcal{X} \to \mathbb{P}^M_T$ by its complete linear series. \textcolor{black}{Note that $M = h^0(L)-1 = h^0(J) + h^0(K_S \otimes J) -1 = N$. Let $\psi: X \to \mathbb{P}^N$ and $\widetilde{\psi}: \widetilde{X} \to \mathbb{P}^N_{\Delta}$ are the morphisms induced by the complete linear series of $L$ and $\widetilde{L}$ respectively. 
Note that $\widetilde{\varphi}$ is obtained by composing $\widetilde{\psi}$ by a projection $\widetilde{p}: \mathbb{P}^N_{\Delta} \to \mathbb{P}^n_{\Delta}$ followed by the linear embedding $k_{\Delta}: \mathbb{P}^n_{\Delta} \hookrightarrow \mathbb{P}^N_{\Delta}$. Let $p: \mathbb{P}^N \to \mathbb{P}^n$ be the restriction of $\widetilde{p}$ to $\mathbb{P}^N$.
Note that we have a cartesian diagram.} 
\color{black}
\[
\begin{tikzcd}
X\arrow[r, "\psi"] \arrow{d} & \mathbb{P}^N \arrow[r, "p"]\arrow{d} & \mathbb{P}^n\arrow[r, "k"]\arrow{d} & \mathbb{P}^N \arrow[r] \arrow{d} & 0\arrow{d}\\
\widetilde{X}\arrow[r, "\widetilde{\psi}"] & \mathbb{P}_{\Delta}^N\arrow[r, "\widetilde{p}"] & \mathbb{P}_{\Delta}^n\arrow[r, "k_{\Delta}"] & \mathbb{P}_{\Delta}^N \arrow{r} & \Delta
\end{tikzcd}\]

Since one can lift the projection $\widetilde{p}$ to a projection $P: \mathbb{P}^M_T \to \mathbb{P}^N_T$ and the embedding $k_{\Delta}$ to an embedding $K: \mathbb{P}^n_T \to \mathbb{P}^N_T$  we get the following diagram. 
\[
\begin{tikzcd}
    X \arrow{r} \arrow[d,"\psi"] & \widetilde{X} \arrow{r} \arrow[d, "\widetilde{\psi}"] & \mathcal{X} \arrow[d, "\Psi"] \\
   \mathbb{P}^N \arrow{r} \arrow[d,"p"] & \mathbb{P}_{\Delta}^N  \arrow{r} \arrow[d, "\widetilde{p}"] & \mathbb{P}_T^N \arrow[d, "P"]  \\
    \mathbb{P}^n \arrow{r} \arrow[d, "k"] & \mathbb{P}_{\Delta}^n  \arrow{r} \arrow[d, "k_{\Delta}"] & \mathbb{P}_T^n \arrow[d, "K"] \\
    \mathbb{P}^N \arrow{r} \arrow{d} & \mathbb{P}_{\Delta}^N  \arrow{r} \arrow{d} & \mathbb{P}_T^N \arrow{d} \\
    0 \arrow{r} &\Delta \arrow{r} & T 
\end{tikzcd} \]
Let $\Phi = K \circ P \circ \Psi$. Then $\Phi$ is a deformation of $\varphi$ which satisfies the conditions of Proposition ~\ref{main}. Hence $\widetilde{S}$ has an embedded smoothing by Proposition ~\ref{main}. \QEDB\par 

\color{black}

\vspace{5pt}

We are left with the embedded K3 carpets on $\mathbb{F}_e$, with $e>2$. We treat them below.

\begin{theorem}\label{smoothing2}
 Let $\widetilde{S}$ be a projective  $K3$ carpet supported on $S = \mathbb{F}_e$ with $e > 2$. Then $\widetilde{S}$ is smoothable.
\end{theorem}

\noindent\textit{Proof.} Let $\widetilde{M}$ be a very ample line bundle on $\widetilde{S}$. Consider the embedding of $\widetilde{i}: \widetilde{S} \hookrightarrow \mathbb{P}^N$ given by the complete linear series $\widetilde{M}$. As in Theorem ~\ref{smoothing1}, the induced embedding $i: S \hookrightarrow \mathbb{P}^N$ satisfies the 
situation of Lemma ~\ref{cohomology} \par
\textcolor{black}{We know by \cite{Sei92}, Theorem $13$, that there exists a deformation $\mathcal{S} \xrightarrow{p'} T$ of $S$ over an irreducible affine curve $T$ such that $\mathcal{S}_0 = S$ and $\mathcal{S}_t = S' \cong \mathbb{F}_e$ with $0 \leq e \leq 2$.} Given the embedding $i$, we have the restriction of Euler sequence to $S$ (here $\mathcal{O}_S(1)=\mathcal{O}_S(aC_0+bf)$) 
\begin{equation}\label{need}
0 \to \mathcal{O}_S \to   \mathcal{O}_S(1)^{\oplus N+1} \to T_{\mathbb{P}^N}|_S \to 0.  
\end{equation}
Since $H^1(\mathcal{O}_S(1)) = H^2(\mathcal{O}_S) = 0$, we have that $H^1(T_{\mathbb{P}^N}|_S) = 0$, consequently the deformation $\mathcal{S} \to T$ can be realized inside $\mathbb{P}_T^N$, i.e. we have the following diagram.
\[
\begin{tikzcd}
\mathcal{S} \arrow[rd, "p'"] \arrow[r, hook] &  \mathbb{P}_T^N \arrow[d] \\
 & T 
\end{tikzcd} \]
By lemma ~\ref{h1}, (2), we have that $H^1(\mathcal{N}_{S/\mathbb{P}^N} \otimes K_S) = 0$. Thus, after possibly shrinking $T$, $p'_*(\mathcal{N}_{\mathcal{S}/\mathbb{P}_T^N} \otimes K_{\mathcal{S}/T})$ is free of rank $h^0(\mathcal{N}_{S/\mathbb{P}^N} \otimes K_S)$. Consider the projective bundle $$\mathbb{P}(p'_*(\mathcal{N}_{\mathcal{S}/\mathbb{P}_T^N} \otimes K_{\mathcal{S}/T})) \xrightarrow{\pi} T.$$ 
By Theorem ~\ref{counting carpets}, (3), we see that there exists a  open set $U \subset \mathbb{P}(p'_*(\mathcal{N}_{\mathcal{S}/\mathbb{P}_T^N} \otimes K_{\mathcal{S}/T}))$ which intersects every fibre of $\pi$ in a nontrivial open set such every $u \in U$ corresponds to a nowhere vanishing section of $H^0(\mathcal{N}_{\mathcal{S}_t/\mathbb{P}^N_t} \otimes K_{\mathcal{S}_t})$. Now $\widetilde{S} \in  H^0(\mathcal{N}_{S/\mathbb{P}^N} \otimes K_S)$ is a nowhere vanishing section and hence is a point in $U$ in the fibre of $0 \in T$. One can then choose a section of $\pi$ through $\widetilde{S}$ to extend the section $\widetilde{S} \in H^0(\mathcal{N}_{S/\mathbb{P}^N} \otimes K_S)$ to a section $\widetilde{S}_T  \in H^0(\mathcal{N}_{\mathcal{S}/\mathbb{P}_T^N} \otimes K_{\mathcal{S}/T}) $. Since $\widetilde{S} \in U$, one can assume after possibly shrinking $T$, that $\widetilde{S}_{T,t}$ is nowhere vanishing for all $t \in T$. Hence fibrewise each section $\widetilde{S}_{T,t} \in H^0(\mathcal{N}_{\mathcal{S}_t/\mathbb{P}^N_t} \otimes K_{\mathcal{S}_t})$ gives a $K3$ carpet on $\mathcal{S}_t$ embedded in $\mathbb{P}_t^N$. Since each $\mathcal{S}_t \cong \widetilde{S}'\cong \mathbb{F}_e$ for $0 \leq e \leq 2$, we have proved the first part of our claim.\par 
The last statement follows since we just showed that there exists an irreducible component of the Hilbert scheme containing $\widetilde{S}$ inside $\mathbb{P}^N$ whose general element is a smooth $K3$ surface. \QEDB\par 
\subsection{Smoothings of embedded K3 carpets on $\mathbb{P}^2$} We treat the embedded K3 carpets on the $d$-uple embeddings of $\mathbb{P}^2$. The following gives the smoothing result for them.
\begin{theorem}\label{smoothingp2}
Let $\widetilde{S}$ be a projective $K3$ carpet on $S = \mathbb{P}^2$. Then $\widetilde{S}$ is smoothable.
\end{theorem}
\noindent\textit{Proof.} The proof is similar to that of Theorem ~\ref{smoothing1} using Theorem ~\ref{countingp2}.\QEDB\par

\vspace{5pt}

The smoothing results and the fact that two K3 surfaces are deformation equivalent shows
\begin{corollary}\label{deformation equivalence}
Suppose $\widetilde{S}_1$ and $\widetilde{S}_2$ be two abstract K3 carpets supported on minimal rational surfaces or $\mathbb{F}_1$. Then $\widetilde{S}_1$ and $\widetilde{S}_2$ are deformation equivalent.
\end{corollary}

\begin{remark}
Note that the abstract split $K3$ carpets on either a Hirzebruch surface $\mathbb{F}_e$ or $\mathbb{P}^2$ are projective by the proof of Theorem ~\ref{counting carpets} since the split carpet corresponds to the zero element $\tau \in H^{1,1}(S)$ and $\int_D \tau = 0$ for any ample divisor $D$ over either $\mathbb{F}_e$ or $\mathbb{P}^2$. Hence they are smoothable.
\end{remark}

The results in \cite{BGG} show that there are no higher dimensional analogues of the above results. They introduce the notion of generalized hyperelliptic varieties that unifies various results on deformations of double covers. We now give that definition. 

\begin{definition}\label{ghpv}
Let $(X, L)$ a smooth polarized variety with $L$ base-point-free. Let the morphism
from $X$ to $\mathbb{P}^N$, induced by $|L|$ be of degree $2$ onto its image $i(Y)$ ($i$ is the embedding of $Y$ in $\mathbb{P}^N$). Let the variety $Y$ be smooth and isomorphic to any of these: (1)  projective space, (2) a hyperquadric, (3) a projective bundle over $\mathbb{P}^1$.
The variety $Y$ is not necessarily embedded by $i$ as a variety of minimal degree in $\mathbb{P}^N$. Then we say
that $(X, L)$ is a generalized hyperelliptic polarized variety.
\end{definition}

Calabi-Yau varieties are higher dimensional analogues of K3 surfaces. In \cite{BGG}, they show that there are no Calabi-Yau double structures on projective bundles and $\mathbb {P}^N$, in sharp contrast with the results on K3 surfaces in this article. This together with results on deformations in \cite{BGG} show that deformations of generalized hyperelliptic Calabi-Yau $n$-fold is again a generalized hyperelliptic $n$-fold.

\section{Hilbert Points of K3 Carpets}
In this section we prove results on the Hilbert point of a projective K3 carpet on the Hirzebruch surfaces and $\mathbb{P}^2$.
 This is in sharp contrast with the results in \cite{GGP} on Hilbert points corresponding to K3
carpets on an Enriques surface, where all the Hilbert points were smooth. 

\subsection{K3 carpets on Hirzebruch surfaces} First we deal with the K3 carpets on $\mathbb{F}_e$. The following result generalizes \cite{GP97}, Theorem 4.1.
\begin{theorem}\label{sm1}
 
Let $\widetilde{S}$ be a projective $K3$ carpet supported on $S = \mathbb{F}_e$ embedded in $\mathbb{P}^N$ by the complete linear series of a very ample line bundle. Then $\widetilde{S}$ is a smooth point of its Hilbert scheme if and only if $0 \leq e \leq 2$. 
\end{theorem}

\noindent\textit{Proof.} We know from Theorem ~\ref{smoothing1} and Theorem ~\ref{smoothing2} that $\widetilde{S}$ admits an embedded smoothing. Hence there exists an irreducible component $H$ of the Hilbert scheme containing $\widetilde{S}$ such that a general element of $H$ is a smooth $K3$ surface $X$. Let $i': X \to \mathbb{P}^N$ denote the embedding. Let $L = \mathcal{O}_X(1):=i'^*\mathcal{O}_{\mathbb{P}^N}(1)$. Then, we have the following exact sequence.
\begin{equation}
0 \to \mathcal{O}_X \to \mathcal{O}_X(1)^{\oplus N+1} \to T_{\mathbb{P}^N}|_X \to 0.
\end{equation}
Using Kodaira vanishing and $h^2(\mathcal{O}_X) = 1$ we have that $h^1(T_{\mathbb{P}^N}|_X) = 1$. Now consider 
\begin{equation}
0 \to T_X \to T_{\mathbb{P}^N}|_X \to \mathcal{N}_{X/\mathbb{P}^N} \to 0.
\end{equation}
Taking cohomology we get that 
\begin{equation}
H^0(\mathcal{N}_{X/\mathbb{P}^N}) \to H^1(T_X) \to H^1(T_{\mathbb{P}^N}|_X) \to H^1(\mathcal{N}_{X/\mathbb{P}^N}) \to H^2(T_X) = 0.
\end{equation}
Note that, since there exist non-projective $K3$ surfaces we have that $H^0(\mathcal{N}_{X/\mathbb{P}^N}) \to H^1(T_X)$ is not surjective. Consequently, $H^1(\mathcal{N}_{X/\mathbb{P}^N}) = 0$. Thus, $X$ is smooth point in its Hilbert scheme. Therefore $\widetilde{S}$ is a smooth point of the Hilbert scheme if and only if $$h^0(\mathcal{N}_{\widetilde{S}/\mathbb{P}^N}) = h^0(\mathcal{N}_{X/\mathbb{P}^N}) =(N+1)^2 + 18.$$
(Note that  $N+1 = h^0(\mathcal{O}_{\widetilde{S}}(1)) = h^0(\mathcal{O}_X(1))$.\par 
In order to compute $h^0(\mathcal{N}_{\widetilde{S}/\mathbb{P}^N})$ we use the following exact sequences; 
\begin{equation}\label{f1}
0 \to \mathcal{N}_{\widetilde{S}/\mathbb{P}^N} \otimes K_S \to \mathcal{N}_{\widetilde{S}/\mathbb{P}^N} \to \mathcal{N}_{\widetilde{S}/\mathbb{P}^N} \otimes \mathcal{O}_S \to 0,
\end{equation}
\begin{equation}\label{f2}
0 \to K_S^* \to \mathcal{N}_{S/\mathbb{P}^N} \to \mathcal{H}om(I_{\widetilde{S}}/I_S^2, \mathcal{O}_S) \to 0,    
\end{equation}
\begin{equation}\label{f3}
0 \to \mathcal{O}_S \to \mathcal{N}_{S/\mathbb{P}^N} \otimes K_S \to \mathcal{H}om(I_{\widetilde{S}}/I_S^2, \mathcal{O}_S) \otimes K_S \to 0,    
\end{equation}
\begin{equation}\label{f4}
0 \to  \mathcal{H}om(I_{\widetilde{S}}/I_S^2, \mathcal{O}_S) \to \mathcal{N}_{\widetilde{S}/\mathbb{P}^N} \otimes \mathcal{O}_S \to K_S^{-2} \to 0,  
\end{equation}
\begin{equation}\label{f5}
0 \to  \mathcal{H}om(I_{\widetilde{S}}/I_S^2, \mathcal{O}_S) \otimes K_S \to \mathcal{N}_{\widetilde{S}/\mathbb{P}^N} \otimes K_S \to K_S^* \to 0.  
\end{equation}
Using (\ref{f2}), we have that $h^1(\mathcal{H}om(I_{\widetilde{S}}/I_S^2, \mathcal{O}_S)) = h^2(\mathcal{H}om(I_{\widetilde{S}}/I_S^2, \mathcal{O}_S)) = 0$, and the following;
$$h^0(\mathcal{H}om(I_{\widetilde{S}}/I_S^2, \mathcal{O}_S)) = h^0(\mathcal{N}_{S/\mathbb{P}^N}) - h^0(K_S^*) + h^1(K_S^*).$$
Using (\ref{f3}), we have that $h^1(\mathcal{H}om(I_{\widetilde{S}}/I_S^2, \mathcal{O}_S) \otimes K_S) = h^2(\mathcal{H}om(I_{\widetilde{S}}/I_S^2, \mathcal{O}_S) \otimes K_S) = 0$, and $$h^0(\mathcal{H}om(I_{\widetilde{S}}/I_S^2, \mathcal{O}_S) \otimes K_S) = h^0(\mathcal{N}_{S/\mathbb{P}^N} \otimes K_S) - 1.$$ 
Using (\ref{f4}), we have that $h^1(\mathcal{N}_{\widetilde{S}/\mathbb{P}^N} \otimes \mathcal{O}_S) = h^1(K_S^{-2}), h^2(\mathcal{N}_{\widetilde{S}/\mathbb{P}^N} \otimes \mathcal{O}_S) = 0$, and the following; $$h^0(\mathcal{N}_{\widetilde{S}/\mathbb{P}^N} \otimes \mathcal{O}_S) = h^0(\mathcal{N}_{S/\mathbb{P}^N}) - h^0(K_S^*) + h^1(K_S^*) + h^0(K_S^{-2}).$$ 
Using (\ref{f5}), we have that $h^1(\mathcal{N}_{\widetilde{S}/\mathbb{P}^N} \otimes K_S) = h^1(K_S^*), h^2(\mathcal{N}_{\widetilde{S}/\mathbb{P}^N} \otimes K_S) = 0$, and the following; $$h^0(\mathcal{N}_{\widetilde{S}/\mathbb{P}^N} \otimes K_S) = h^0(\mathcal{N}_{S/\mathbb{P}^N} \otimes K_S) - 1 + h^0(K_S^*).$$ 
Using (\ref{f1}), we have that 
$h^0(\mathcal{N}_{\widetilde{S}/\mathbb{P}^N}) - h^1(\mathcal{N}_{\widetilde{S}/\mathbb{P}^N}) = h^0(\mathcal{N}_{S/\mathbb{P}^N}) + h^0(\mathcal{N}_{S/\mathbb{P}^N} \otimes K_S) + 24$. It follows from Lemma ~\ref{h1} that $h^0(\mathcal{N}_{S/\mathbb{P}^N} \otimes K_S)=(N+1)h^0(\mathcal{O}_S((aC_0+bf)\otimes K_S))+1$. 
\color{black}
Using the following exact sequence $$0\to T_S\to T_{\mathbb{P}^N}|_S\to\mathcal{N}_{S/\mathbb{P}^N}\to 0$$ and the restriction of the Euler sequence on $Y$ (see ~\eqref{need}), we obtain that $$ h^0(\mathcal{N}_{S/\mathbb{P}^N}) = h^0(T_{\mathbb{P}^N}|_S)-(h^0(T_S)-h^1(T_S))=(N+1)(h^0(\mathcal{O}_S(aC_0+bf)) -7,$$
\color{black}
and hence $h^0(\mathcal{N}_{\widetilde{S}/\mathbb{P}^N}) - h^1(\mathcal{N}_{\widetilde{S}/\mathbb{P}^N}) = (N+1)^2 + 18$ (notice $h^0(\mathcal{O}_S(aC_0+bf)+h^0(\mathcal{O}_S((aC_0+bf)\otimes K_S)) = N+1$). Thus, $\widetilde{S}$ is smooth if and only if $h^1(\mathcal{N}_{\widetilde{S}/\mathbb{P}^N}) = 0$.\par 

Using (\ref{f1}), we get the following exact sequence;
\begin{equation}
H^1(\mathcal{N}_{\widetilde{S}/\mathbb{P}^N} \otimes K_S) \cong H^1(K_S^*) \to H^1(\mathcal{N}_{\widetilde{S}/\mathbb{P}^N}) \to H^1(\mathcal{N}_{\widetilde{S}/\mathbb{P}^N} \otimes \mathcal{O}_S) \cong H^1(K_S^{-2}) \to 0.
\end{equation}
For $0 \leq e \leq 2$ both the flanking terms vanish but for $e > 2$, $h^1(K_S^{-2}) \neq 0$. Thus, $\widetilde{S}$ is smooth if and only if $0 \leq e \leq 2$. \QEDB\par

\subsection{K3 carpets on $\mathbb{P}^2$} We conclude the article by the following result. 
\begin{theorem}\label{sm2}
 
Let $\widetilde{S}$ be a projective $K3$ carpet supported on $S=\mathbb{P}^2$ and embedded inside $\mathbb{P}^N$ by the complete linear series of a very ample line bundle. Then $\widetilde{S}$ is a smooth point of its Hilbert scheme. 
\end{theorem}
\noindent\textit{Proof.} As in the previous theorem we need to show $$h^0(\mathcal{N}_{\widetilde{S}/\mathbb{P}^N}) = (N+1)^2 + 18.$$ 
Using (\ref{f2}) and Lemma ~\ref{p2}, $(1)$ we have that $h^1(\mathcal{H}om(I_{\widetilde{S}}/I_S^2, \mathcal{O}_S)) = h^2(\mathcal{H}om(I_{\widetilde{S}}/I_S^2, \mathcal{O}_S)) = 0$, and the following;  $$h^0(\mathcal{H}om(I_{\widetilde{S}}/I_S^2, \mathcal{O}_S)) = (N+1)h^0(\mathcal{O}_{\mathbb{P}^2}(d)) - 3h^0(\mathcal{O}_{\mathbb{P}^2}(1)) - h^0(\mathcal{O}_{\mathbb{P}^2}(3)).$$ 
Using (\ref{f3}) and ~\eqref{3.1}, we have that $h^1(\mathcal{H}om(I_{\widetilde{S}}/I_S^2, \mathcal{O}_S) \otimes K_S) = h^2(\mathcal{H}om(I_{\widetilde{S}}/I_S^2, \mathcal{O}_S) \otimes K_S) = 0$, and $$h^0(\mathcal{H}om(I_{\widetilde{S}}/I_S^2, \mathcal{O}_S) \otimes K_S) = (N+1)h^0(\mathcal{O}_{\mathbb{P}^2}(d-3)) - 1.$$ 
Using (\ref{f4}), we have that $h^1(\mathcal{N}_{\widetilde{S}/\mathbb{P}^N} \otimes \mathcal{O}_S) = h^2(\mathcal{N}_{\widetilde{S}/\mathbb{P}^N} \otimes \mathcal{O}_S) = 0$, and the following; $$h^0(\mathcal{N}_{\widetilde{S}/\mathbb{P}^N} \otimes \mathcal{O}_S) = (N+1)h^0(\mathcal{O}_{\mathbb{P}^2}(d)) - 3h^0(\mathcal{O}_{\mathbb{P}^2}(1)) - h^0(\mathcal{O}_{\mathbb{P}^2}(3))+ h^0(\mathcal{O}_{\mathbb{P}^2}(6)). 
$$ 
Using (\ref{f5}), we have that $h^1(\mathcal{N}_{\widetilde{S}/\mathbb{P}^N} \otimes K_S) = h^1(\mathcal{O}_{\mathbb{P}^2}(3)) = 0, h^2(\mathcal{N}_{\widetilde{S}/\mathbb{P}^N} \otimes K_S) = h^2(\mathcal{O}_{\mathbb{P}^2}(3)) = 0$, and 
$$h^0(\mathcal{N}_{\widetilde{S}/\mathbb{P}^N} \otimes K_S) = (N+1)h^0(\mathcal{O}_{\mathbb{P}^2}(d-3)) - 1 + h^0(\mathcal{O}_{\mathbb{P}^2}(3)).$$ 
Using (\ref{f1}), we have that 
$h^0(\mathcal{N}_{\widetilde{S}/\mathbb{P}^N}) = h^0(\mathcal{N}_{\widetilde{S}/\mathbb{P}^N} \otimes \mathcal{O}_S) + h^0(\mathcal{N}_{\widetilde{S}/\mathbb{P}^N} \otimes K_S) $. It follows 
that $$h^0(\mathcal{N}_{\widetilde{S}/\mathbb{P}^N}) = (N+1)(h^0(\mathcal{O}_{\mathbb{P}^2}(d-3)) + h^0(\mathcal{O}_{\mathbb{P}^2}(d))) + 18.$$ The assertion follows since $(N+1) = h^0(\mathcal{O}_{\mathbb{P}^2}(d-3)) + h^0(\mathcal{O}_{\mathbb{P}^2}(d))$.\QEDB

\bibliographystyle{plain}

\end{document}